\documentclass[12pt]{article}
\pdfoutput=1
\usepackage{latexsym, graphicx,cite,hyphenat}
\usepackage{amsmath}
\usepackage{amssymb}
\usepackage{amsthm}
\allowdisplaybreaks

\usepackage[top = 1. in,bottom = 1. in, left = 1 in, right=1 in]{geometry}

 \newcommand{\arXiv}[1]{\href{http://www.arXiv.org/abs/#1}{#1}}
\usepackage[colorlinks=true, linkcolor=blue, bookmarks=true]{hyperref}

\makeatletter
\renewcommand\section{\@startsection {section}{1}{\z@}%
                  {-3.5ex \@plus -1ex \@minus -.2ex}
                  {2.3ex \@plus.2ex}%
                  {\normalfont\large\bfseries}}
\renewcommand\subsection{\@startsection{subsection}{2}{\z@}%
                   {-3.25ex\@plus -1ex \@minus -.2ex}%
                   {1.5ex \@plus .2ex}%
                   {\normalfont\bfseries}}
\makeatother

\newcommand{\beq}{\begin{equation}}
\newcommand{\eeq}{\end{equation}}
\newcommand{\ber}{\begin{array}}
\newcommand{\eer}{\end{array}}
\newcommand{\del}{\partial}

\newcommand{\ssty}{\scriptstyle}

\newcommand{\s}{y} 
\newcommand{\te}{\theta}

\newcommand{\eps}{\varepsilon}

\newcommand{\ena}{\end{eqnarray}}
\newcommand{\beqa}{\begin{eqnarray}}
\newcommand{\eeqa}{\end{eqnarray}}
\newcommand{\bea}{\begin{eqnarray}}
\newcommand{\eea}{\end{eqnarray}}

\theoremstyle{remark}
\newtheorem{rem}{Remark}

\renewcommand{\Re}{\operatorname{Re}}
\renewcommand{\Im}{\operatorname{Im}}

\newcommand{\minPlusOne}[1]{
  \min\!\left(#1\right) + 1
}

\begin{document}
\begin{titlepage}
\begin{flushright}
\phantom{arXiv:yymm.nnnn}
\end{flushright}
\vspace{-2cm}
\begin{center}
{\LARGE\bf \nohyphens{Conformal flow on S$^3$ and}
\vspace{4mm}
\\weak field integrability in AdS$_4$}  \\
\vskip 15mm
{\large Piotr Bizo\'n$^{a,d}$, Ben Craps$^{b}$, Oleg Evnin$^{c,b}$,\vspace{1mm}\\
Dominika Hunik$^a$, Vincent Luyten$^{b}$, Maciej Maliborski$^{d}$}
\vskip 7mm
{\em $^a$ Institute of Physics, Jagiellonian University, Krak\'ow, Poland}
\vskip 3mm
{\em $^b$ Theoretische Natuurkunde, Vrije Universiteit Brussel and\\
The International Solvay Institutes, Brussels, Belgium}
\vskip 3mm
{\em $^c$ Department of Physics, Faculty of Science, Chulalongkorn University,
Bangkok, Thailand}
\vskip 3mm
{\em $^d$ Max-Planck-Institut f\"ur Gravitationsphysik, Albert-Einstein-Institut, Golm, Germany}
\vskip 7mm
{\small\noindent {\tt bizon@th.if.uj.edu.pl, Ben.Craps@vub.ac.be, oleg.evnin@gmail.com,\\ dominika.hunik@uj.edu.pl, Vincent.Luyten@vub.ac.be, maciej.maliborski@aei.mpg.de}}
\vskip 10mm
\end{center}
\vspace{1cm}
\begin{center}
{\bf ABSTRACT}\vspace{3mm}
\end{center}
We consider the conformally invariant cubic wave equation on the Einstein cylinder $\mathbb{R} \times \mathbb{S}^3$ for small rotationally symmetric initial data. This simple equation captures many key challenges of nonlinear wave dynamics in confining geometries, while a conformal transformation relates it to a self-interacting conformally coupled scalar in four-dimensional anti-de Sitter spacetime (AdS$_4$) and connects it to various questions of AdS stability. We construct an effective infinite-dimensional time-averaged dynamical system accurately approximating
the original equation in the weak field regime. It turns out that this effective system, which we call the \emph{conformal flow}, exhibits some remarkable features,
such as low-dimensional invariant subspaces, a wealth of stationary
states (for which energy does not flow between the modes), as well as solutions with nontrivial exactly periodic energy
flows. Based on these observations and close parallels to the cubic Szeg\H{o} equation, which was shown by G\'erard and Grellier to be Lax-integrable, it is tempting to conjecture that the conformal flow and the corresponding weak field dynamics in AdS$_4$ are integrable as well.

\vfill

\end{titlepage}


\section{Introduction}

Propagation of nonlinear waves in confining geometries presents significant challenges
because the key mechanism stabilizing the evolution of waves on unbounded domains, which is the dispersion of energy by radiation, is missing in
confined settings. Consequently, an arbitrarily small perturbation of a ground state can lead to complicated long-time behavior. The
central physical problem in this context is that of energy transfer,
namely, understanding how the energy injected into the system gets distributed over the degrees
of freedom in the course of its evolution. This problem has been studied in
the physics literature from a statistical viewpoint under the name
of wave (weak) turbulence \cite{Zakharov,Nazarenko}, but only recently emerged as
an active research topic in dispersive PDE theory. In this deterministic mathematical approach, the energy cascades from low to high modes are quantified by the growth of higher Sobolev norms of solutions and the main question is whether these norms can become unbounded in finite or infinite time. The past few years
have witnessed significant progress in the
understanding of this issue in the context of nonlinear Sch\"odinger equations (and their variations) on tori \cite{KSS, H, GK}. However, almost nothing is known, as far as we can tell, about manifestations
of these phenomena in other interesting evolution equations of mathematical physics.

Our attention in this article will be focused on the weak field dynamics for the conformally invariant cubic wave equation on the three-sphere. This very natural geometrical model describes relativistic scalar waves propagating on a compact manifold and interacting with the background geometry and among themselves. The spherical geometry prevents the waves from dispersing to infinity and thus ensures that the nonlinear self-interactions remain important for all times, inducing complicated energy transfer patterns. Another important circumstance is that all normal mode frequencies of the corresponding linearized theory are integer in appropriate units. Such a fully resonant spectrum ensures that nonlinearities can produce significant effects over long times for arbitrarily small perturbation amplitudes, generating highly nontrivial weak field dynamics.

The weak field dynamics of the system we are considering is closely related to investigations of nonlinear stability of anti-de Sitter (AdS) spacetime initiated in \cite{BR} (see \cite{Bizon:2013gxa, review} for brief reviews and further references). A conformal transformation relates our considerations to the dynamics of the cubic Klein-Gordon equation with mass $m^2=-2$ (conformally coupled self-interacting probe scalar) on the four-dimensional AdS space. In contrast to the bulk of AdS stability research, we do not consider the backreaction of perturbations on the AdS geometry. However, a number of features of the weak field dynamics remain unaffected by our simplification. (The weak field dynamics of a massless self-interacting probe scalar in AdS was introduced as a toy model for the AdS stability problem in \cite{BKS1}.)

Our main tool for analyzing the weak field dynamics, given a fully resonant spectrum of the linearized system, is the time averaging method that goes back to Bogoliubov and Krylov and is described in many books on perturbation theory (e.g., \cite{murdock}). For a contemporary treatment adapted to studies of nonlinear PDEs, see \cite{KM}. Applied to the PDE we are considering, the time averaging method produces a simplified infinite-dimensional dynamical system, which we shall refer to as the {\em conformal flow}. By standard theorems underlying time averaging, this effective system (in PDE theory sometimes called the \emph{resonant} system) accurately approximates the dynamics of the original PDE for small amplitude fields of order $O(\eps)$ on long time-scales of order $O(1/\eps^2)$, and in particular it accurately captures the energy transfer patterns on such time scales. (Time averaging and related techniques were introduced to the AdS stability problem closely related to our current studies in \cite{BBGLL, CEV1,CEV2}.) The conformal flow displays a number of highly special features not apparently present in the original PDE and suggestive of integrability. This includes additional conserved quantities, low-dimensional invariant subspaces, a variety of stationary states without any energy transfer, as well as solutions with exactly periodic energy transfer patterns.

Our original PDE can be viewed as an infinite system of oscillators with integer frequencies and a quartic potential. We do not see any indications that it might be integrable; only the time-averaged system describing its weak field dynamics displays signs of integrability. While extensive searches for quartic integrable systems (i.e., integrable mechanical analogs of our original PDE) have been conducted in the past (see \cite{LS} for a review and references to earlier work and \cite{BC} for some later contributions), we are not aware of any broad scans for systems displaying {\em weak field integrability} of the type suggested by our results.

The conformal flow introduced in this paper bears a strong resemblance to the cubic Szeg\H{o} equation that has been designed and studied in a series of papers by G\'erard and Grellier, of which \cite{GG, GGint, GG2, GG3} are particularly relevant for our purposes. This equation also emerges  as the resonant approximation of a nonlocal transport equation, called the cubic half-wave equation \cite{GGint}, which is a special case of one-dimensional models of wave turbulence introduced  in \cite{MLT}.
In\cite{GG, GGint, GG2, GG3} and related publications, a number of remarkable properties were established for  the cubic Szeg\H{o} equation, in particular, a Lax pair structure, existence of finite-dimensional invariant subspaces, as well as some weak turbulence phenomena. The bulk of our effort will be directed at recovering analogs of a subset of these results for the conformal flow, a system with a natural geometric origin. We furthermore believe that the analogy goes deeper than what we are able to explicitly demonstrate in this article.

There are also close similarities between the conformal flow and   resonant approximations for the cubic Schr\"odinger equation, either  with harmonic trapping  \cite{HT} or in the infinite volume limit on the 2-dimensional torus \cite{FGH, GHT}. Another system which shares some properties with the conformal flow and the cubic Szeg\H{o} equation is  the Lowest Landau Level evolution equation which appears  in the studies of rapidly rotating Bose-Einstein condensates (see \cite{GT} and references therein).
 We furthermore mention the curious case of the Charney-Hasegawa-Mima (CHM) equation explored in \cite{KM}. Also in that case, effective integrable dynamics emerges in the time-averaged approximation (which the authors call the `effective equation'). However, a substantial difference is that the linearized spectrum of the CHM equation is not fully resonant and splits into disconnected resonant clusters. In such a situation, the resulting effective time-averaged system literally separates into independent finite-dimensional integrable subsystems. This simplification does not happen in our case because our linearized spectrum is fully resonant.

We comment on how our work is related to extensive investigations of the last decade in the area of `AdS/CFT Integrability' (see \cite{AdSreview} for a review). One of the central ingredients of that line of research is integrability of sigma models (describing string worldsheets) on target spaces involving AdS factors. While we are also presently talking about integrability emerging due to special features of the AdS geometry, there are also essential differences. We are considering a field system evolving in AdS, rather than a sigma model (whose fields are maps from a two-dimensional surface into AdS). We are furthermore not talking about exact integrability of our AdS field system, but rather about the integrability of the effective weak field theory (conformal flow) emerging from it within the standard time-averaged approximation.

Our exposition is organized as follows. In section~2, we formulate our scalar field dynamics on a three-sphere, describe its connection to AdS stability problems and apply time averaging to derive the corresponding effective system, which we name the conformal flow. In section~3, we give a very elementary pragmatic summary of some properties of the cubic Szeg\H{o} equation relevant in our context, aimed mostly at the physics audience. In section~4, we construct a three-dimensional invariant subspace of the conformal flow and discuss some complex plane properties of the generating functions that encode the conformal flow amplitudes. In section~5, we analyze stationary states of the conformal flow, in which the initial conditions are adjusted to ensure that no energy transfer occurs. In section~6, we exhibit explicit solutions with periodic energy flows on the three-dimensional invariant subspace of the conformal flow and discuss their properties. Finally, in section~7, we give a summary of what we have practically demonstrated and an outlook in more ambitious directions, keeping the likelihood of Lax-integrability in mind.


\section{Conformally coupled scalar field on the Einstein cylinder\\
and in AdS$_4$ and its weak field dynamics}

As a simple model of confining geometry we consider the Einstein cylinder, which is the globally hyperbolic spacetime $(\mathcal{M},g)$ with topology $\mathbb{R}\times \mathbb{S}^3$ and metric
\begin{equation}\label{metric}
   g= -dt^2 + r^2 \left(dx^2 +\sin^2{\!x}\, d\omega^2\right)\,,
  \end{equation}
where $x\in [0,\pi]$, $d\omega^2$ is the round metric on the unit 2-sphere and $r$ is the radius of $\mathbb{S}^3$. This spacetime has a constant scalar curvature $R(g)=6/r^2$.

As a model of nonlinear dynamics on $\mathcal{M}$, we consider the semilinear wave equation for the real scalar field $\phi:\mathcal{M}\mapsto \mathbb{R}$
 \begin{equation}\label{eq-conf}
  \left(\square_g -\frac{1}{6} R(g)\right) \phi -\lambda \phi^3 =0\,,
  \end{equation}
  where $\square_g:=g^{\mu\nu}\nabla_{\mu}\nabla_{\nu}$ is the wave operator associated with $g$ and $\lambda$ is a constant. For concreteness, we assume that $\lambda>0$ (which corresponds to a defocusing nonlinearity); however, in the small data regime that we focus upon here, all the results below hold true in the focusing case $\lambda<0$ as well. Due to the identity
  \begin{equation}\label{conf-inv}
   \left(\square_{\Omega^2 g} - \frac{1}{6} R(\Omega^2 g)\right) (\Omega^{-1} \phi)=
   \Omega^{-3} \left(\square_{g} - \frac{1}{6} R(g)\right) \phi\,,
  \end{equation}
  equation \eqref{eq-conf} is conformally invariant.
  After rescaling $t\rightarrow t/r$ and $\phi \rightarrow r \sqrt{\lambda}\, \phi$, it takes the dimensionless form
  \begin{equation}\label{eq-conf2}
  \phi_{tt}-\Delta_{\mathbb{S}^3} \phi + \phi + \phi^3 =0\,,
  \end{equation}
   which can be interpreted as the cubic Klein-Gordon equation on the unit $\mathbb{S}^3$ with a unit mass. For simplicity, we shall restrict our analysis to rotationally symmetric fields, which depend only on $(t,x)$.
   Substituting $\phi(t,x)=v(t,x)/\sin(x)$ into \eqref{eq-conf2} we obtain a nonlinear string equation
\begin{equation}\label{eq}
v_{tt} - v_{xx}+\frac{v^3}{\sin^2{x}}=0\,,
\end{equation}
with Dirichlet boundary conditions $v(t,0)=v(t,\pi)=0$, which are enforced by regularity of $\phi(t,x)$ on $\mathcal{M}$.
\begin{rem}
We could have arrived at the same equation considering the conformally coupled self-interacting scalar (\ref{eq-conf}) in the 4-dimensional anti-de Sitter (AdS$_4$) spacetime with the metric
\beq
\tilde g=\frac{1}{\cos^2 x}\left(-dt^2+dx^2+\sin^2{\!x}\, d\omega^2\right),
\label{adsmetric}
\eeq
where $x$ varies between 0 and $\pi/2$. Indeed, $g = \cos^2{x}\, \tilde g$, hence it follows from \eqref{conf-inv} that the field redefinition $\tilde\phi(t,x)=\phi(t,x)\cos x$ converts equation (\ref{eq-conf}) from the Einstein cylinder to AdS. In order to define the evolution, one has to impose a boundary condition at $x=\pi/2$. Here, we shall impose the Dirichlet condition $\phi(\pi/2)=0$ on the equator, corresponding to the reflecting boundary conditions most commonly used in the AdS research context. This model is just a subsector of the model on the full Einstein cylinder, which can be implemented by imposing the reflection symmetry $\phi(t, x)=-\phi(t,\pi-x)$. We emphasize this connection because problems involving nonlinear dynamics of small AdS perturbations have received a significant amount of attention in recent years, and our objectives here have much in common with that body of work.
\end{rem}

 Our goal is to understand the evolution of small smooth initial data $v(0,x), \del_t v(0,x)$.
Decomposing the solution into a Fourier series,
\beq
v(t,x)=\sum\limits_{n=0}^{\infty} c_n(t) \sin(n+1) x,
\eeq
we get from \eqref{eq} an infinite system of coupled oscillators
\begin{equation}\label{fourier}
 \frac{d^2 c_n}{dt^2} +(n+1)^2 c_n = -\sum\limits_{jkl} S_{jkln}\, c_j c_k c_l\,,
\end{equation}
with the interaction coefficients
\begin{equation}\label{C}
 S_{jkln}=\frac{2}{\pi}\,\int_0^{\pi} \frac{dx}{\sin^2{x}} \sin(j+1)x\,\sin(k+1)x\, \sin(l+1)x\, \sin(n+1)x \,.
\end{equation}
To factor out fast linear oscillations in \eqref{fourier}, we change the variables using variation of constants
\begin{eqnarray}\label{voc}
 c_n &=&\beta_n e^{i (n+1) t} +\bar \beta_n e^{-i (n+1) t},\\
 \frac{d c_n}{dt} &=& i (n+1) \left(\beta_n e^{i(n+1) t} -\bar \beta_n e^{-i (n+1) t}\right)\,.
\end{eqnarray}
This transforms the system \eqref{fourier} into
\begin{equation}\label{eqs}
 2 i (n+1) \frac{d \beta_n}{dt} = -\sum\limits_{jkl} S_{jkln} \,
 c_j c_k c_l \, e^{-i (n+1) t}\,,
 \end{equation}
 where each $c_j$ in the sum is given by (\ref{voc}). Thus, each term in the sum
 has a factor $e^{-i\Omega t}$, where $\Omega=(n+1)\pm (j+1)\pm (k+1) \pm (l+1)$, with all the three plus-minus signs independent. The terms with $\Omega=0$ correspond to resonant interactions, while those with $\Omega\neq 0$ are non-resonant.

 Passing to slow time $\tau=\varepsilon^2 t$ and rescaling $\beta_n(t)=\varepsilon \alpha_n(\tau)$, we see that
for $\varepsilon$ going to zero the non-resonant terms  $\propto e^{-i\Omega \tau/\varepsilon^2}$ are highly oscillatory, and therefore expected to be negligible (in a sense we are about to specify). Keeping only the resonant terms in \eqref{eqs} (which is equivalent to time-averaging), we obtain an infinite autonomous dynamical system
\begin{equation}\label{rs}
 2 i (n+1) \frac{d \alpha_n}{d\tau} = -3 \sum\limits_{j k l} S_{jkln} \,\alpha_j \alpha_k \bar \alpha_l \,,
\end{equation}
where the summation runs over the set of indices $\{j,k,l\}$ for which $\Omega=0$.  Evaluating the integrals for the $S$-coefficients, one finds that this set reduces to $\{jkl\, |\, j+k-l=n\}$ and for such resonant combinations of indices,
\begin{equation}\label{S}
 S_{jkln} = \minPlusOne{j,k,l,n}\,.
\end{equation}
Note that terms with $n=j+k+l+2$ could have been present in principle, but the corresponding $S$-coefficients vanish, which can be verified by direct calculation. This is directly parallel to the selection rules that have been extensively discussed in the AdS stability literature \cite{CEV1,CEV2,Yang,EN}.

 It is part of the standard lore in nonlinear perturbation theory that solutions to \eqref{eqs} starting from small initial data of size $\mathcal{O}(\varepsilon)$ are well approximated by solutions of \eqref{rs} on timescales of order $\mathcal{O}(\varepsilon^{-2})$. More precisely, if $\beta_n(t)$ and $\alpha_n(\tau)$ are solutions of \eqref{eqs} and \eqref{rs}, respectively, and $\beta_n(0)=\varepsilon \alpha_n(0)$ for each $n$, then $|\beta_n(t)-\varepsilon \alpha_n(\tau(t))| \lesssim \mathcal{O}(\varepsilon^2)$ for $t \lesssim \mathcal{O}(\varepsilon^{-2})$.
In other words, on this timescale the dynamics of solutions of equation \eqref{eq} is dominated by resonant interactions. Straightforward proofs (normally phrased for systems with a finite number of degrees of freedom) can be found in textbooks on nonlinear perturbation theory, e.g., \cite{murdock}.

For convenience, we shall henceforth rescale $\tau$ to remove the numerical factors 2 and $-3$ from (\ref{rs}). Renaming so rescaled $\tau$ back to $t$ and using $\cdot=d/dt$, we finally arrive at the following system, which we call the {\em conformal flow}
\beq
i(n+1)\dot \alpha_n = \sum\limits_{j=0}^\infty \sum_{k=0}^{n+j} [\minPlusOne{n,j,k,n+j-k}]\bar\alpha_j \alpha_k \alpha_{n+j-k}\,.
\label{flow}
\eeq

\begin{rem}
Note that the AdS formulation described in \emph{Remark 1} simply corresponds to setting all even-numbered modes to 0, while only keeping odd-numbered modes $\alpha_{2m+1}$ with $m$ running from 0 to infinity. It is easy to see that if one implements this constraint in (\ref{flow}) and expresses everything through the new mode counting index $m$, one gets back equation (\ref{flow}). Thus, the conformal flow is equally applicable to the dynamics in  $\mathbb{R} \times \mathbb{S}^3$ and in AdS$_4$.
\end{rem}

The conformal flow \eqref{flow} is Hamiltonian with
\begin{equation}
 H= \sum\limits_{n=0}^{\infty}\sum\limits_{j=0}^{\infty} \sum\limits_{k=0}^{n+j}[\minPlusOne{n,j,k,n+j-k}]\bar\alpha_n \bar \alpha_j \alpha_k \alpha_{n+j-k}
\label{Hconf}
\end{equation}
and symplectic form $\sum_n 2i(n+1)\,d\bar\alpha_n\wedge d\alpha_n$:
\begin{equation}\label{ham-eq}
i (n+1) \dot \alpha_n = \frac{1}{2}\,\frac{\partial H}{\partial \bar \alpha_n}\,.
\end{equation}
 It enjoys the following one-parameter groups of symmetries (where $\lambda$, $\theta$ are real parameters):
 \begin{eqnarray}
 \mbox{Scaling:} \,&& \alpha_n(t) \rightarrow \lambda \alpha_n(\lambda^2 t),\label{symmscale}\\
 \mbox{Global phase shift:} \,&& \alpha_n(t) \rightarrow e^{i\theta} \alpha_n(t),\label{symmshift1}\\
 \mbox{Mode-dependent phase shift:} \, &&\alpha_n(t) \rightarrow e^{i n \theta} \alpha_n(t).\label{symmshift2}
 \end{eqnarray}
The latter two symmetries respect the Hamiltonian structure and give rise to two conserved quantities (in addition to the Hamiltonian itself):
\begin{eqnarray}
\mbox{`Charge:'} \,&& Q =\sum\limits_{n=0}^{\infty} (n+1) |\alpha_n|^2,\label{charge}\\
\mbox{`Linear energy:'} \,&& E =\sum\limits_{n=0}^{\infty} (n+1)^2 |\alpha_n|^2\,. \label{lenergy}
\end{eqnarray}
We emphasize that these two conservation laws have no obvious counterparts in the original wave equation (\ref{eq-conf2}). Their prototypes were presented in the AdS stability context in \cite{BKS1,CEV2}.


\section{An elementary introduction to the cubic Szeg\H{o} equation}\label{szegosection}

The main content of our current treatment will be in presenting a range of remarkably simple dynamical behaviors of the conformal flow system (\ref{flow}). In anticipation of these results, it is useful to contemplate for a moment the apparent special features of equation (\ref{flow}) itself.

One can trivially rewrite (\ref{flow}) as
\beq
i\dot \alpha_n = \sum\limits_{j=0}^\infty \sum_{k=0}^{n+j} \frac{\minPlusOne{n,j,k,n+j-k}}{n+1}\bar\alpha_j \alpha_k \alpha_{n+j-k}.
\label{flowdenom}
\eeq
One notable property of this representation is that many coefficients of the trilinear form on the right hand side are simply 1. Indeed, this happens for the contribution of any triplet of modes with all frequencies higher than the recepient mode $n$. It turns out useful to consider the following simpler equation, in which {\em all} numerical coefficients are set to 1,
\beq
i\dot\alpha_n = \sum\limits_{j=0}^\infty \sum_{k=0}^{n+j} \bar\alpha_j \alpha_k \alpha_{n+j-k}.
\label{GGalpha}
\eeq
We shall see that, while the systems \eqref{flowdenom} and \eqref{GGalpha} are distinct, their dynamics is qualitatively similar, with a number of parallels between the emerging algebraic structures.

Actually, the simplified system \eqref{GGalpha} is the Fourier representation of a simple paradifferential equation that has been studied before. To show this, let $u(t,e^{i\theta})$ be a complex function on the circle whose nonnegative Fourier coefficients are given by $\alpha_n(t)$ and negative Fourier coefficients vanish, i.e.
\begin{equation}\label{ucircle}
  u(t, e^{i\theta})=\sum\limits_{n=0}^{\infty} \alpha_n(t) e^{i n \theta},
\end{equation}
and assume that $\sum\limits_{n=0}^{\infty} |\alpha_n|^2 < \infty$. Mathematically, this means that $u$ belongs to the Hardy  space  on the  circle $L_{+}^2(S^1)\subset L^2(S^1)$ \cite{G}.
Furthermore, we define an orthogonal projector $\Pi: L^2(S^1)\mapsto L_{+}^2(S^1)$ (called the Szeg\H{o} projector) that filters out negative frequencies
\beq
\Pi\left(\sum_{n=-\infty}^{\infty} \alpha_n e^{in\theta}\right):=\sum_{n=0}^{\infty} \alpha_n e^{in\theta}.
\eeq
Then, it is easy to see that the system \eqref{GGalpha} is equivalent to
\beq
i\del_t u =\Pi(|u|^2u).
\label{GG}
\eeq
This equation  was introduced by G\'erard and Grellier in \cite{GG} under the name of the cubic Szeg\H{o} equation. They showed that this equation has a remarkably deep structure, including a Lax pair, finite-dimensional invariant subspaces and weakly turbulent solutions. We shall demonstrate below that the conformal flow exhibits at least some of this structure (in a modified form), and conjecture that the similarity goes even further. In view of these parallels, we shall present below a very elementary summary of the features of the cubic Szeg\H{o} equation that are important to us.  Readers interested in more details  are referred to the original publications \cite{GG,GGint,GG2,GG3}.

The cubic Szeg\H{o} equation is Hamiltonian with
\begin{equation}\label{H_GG}
 H= \sum\limits_{n=0}^{\infty}\sum\limits_{j=0}^{\infty} \sum\limits_{k=0}^{n+j} \bar\alpha_n \bar \alpha_j \alpha_k \alpha_{n+j-k}\,=\frac{1}{2\pi} \int\limits_0^{2\pi} |u|^4 d\theta
\end{equation}
and symplectic form $\sum_n 2i\,d\bar\alpha_n\wedge d\alpha_n$,
and has the same symmetries
 (\ref{symmscale}-\ref{symmshift2}) as (\ref{flow}). Hence, it possesses two conserved quantities analogous to \eqref{charge} and \eqref{lenergy},
 in addition to the Hamiltonian itself:
\begin{eqnarray}
&& M =\sum\limits_{n=0}^{\infty} |\alpha_n|^2=\frac{1}{2\pi} \int\limits_0^{2\pi} |u|^2 d\theta,\\
&& P =\sum\limits_{n=0}^{\infty} n |\alpha_n|^2=\frac{1}{2\pi} \int\limits_0^{2\pi} (-i \partial_{\theta} u) \bar u \,d\theta.
\end{eqnarray}
Borrowing terminology from studies of the nonlinear Schr\"odinger equation, we shall refer to these charges as `mass' and `momentum,' respectively.

Note that $M+P=\lVert u \rVert_{H_{1/2}}^2$, where $H_s\in L_+^2$ is the Sobolev space on the circle equipped with the norm
\begin{equation}\label{sobolev}
\lVert u \rVert_{H_{s}}^2=\sum\limits_{n=0}^{\infty} (n+1)^{2s} |\alpha_n(t)|^2.
\end{equation}
This fact can be used to control the growth in time  of higher Sobolev norms and thereby prove that the Cauchy problem for the Szeg\H{o} equation is globally well posed for smooth initial data in $H_s$ with $s>1/2$ \cite{GG}.

Any function $u(e^{i\theta})$ of the form (\ref{ucircle}) in $L_{+}^2(S^1)$ can be identified with a holomorphic function $u(z)$ inside the unit disk $|z|<1$ in the complex plane \cite{G}, hence we can write
\beq
u(t,z)=\sum_{n=0}^\infty
\alpha_n(t) z^n.
\label{defu}
\eeq
 In terms of $u(t,z)$, the cubic Szeg\H{o} equation \eqref{GG} takes the form
\beq
i\del_t u(t,z)=\frac1{2\pi i}\oint\limits_{|w|=1} \frac{dw}{w-z}\, \tilde u(t,w) u(t,w)^2,
\label{GGgen}
\eeq
where
\beq
\tilde u(t,z) = \sum_{n=0}^\infty \bar\alpha_n z^{-n}.
\label{defutilde}
\eeq
General properties of the Cauchy-type integral appearing in (\ref{GGgen}) suggest the existence of meromorphic solutions for $u$ with a finite number of isolated poles moving in the complex plane outside the unit disk. Indeed, it was shown in \cite{GG} that the cubic Szeg\H{o} equation admits finite-dimensional invariant subspaces given by rational functions of arbitrarily high degree. Moreover, due to the Lax-integrability (and, consequently, infinitely many conserved quantities), the Szeg\H{o} flow restricted to such subspaces is completely integrable in the sense of Liouville. The simplest invariant manifold with nontrivial dynamics is parametrized by the following single pole ansatz
\begin{equation}\label{ansatz-gg}
 u(t,z)=\frac{b(t)+a(t)z}{1-p(t)z},\qquad |p|<1.
\end{equation}
We shall now demonstrate explicitly by elementary means that the subspace given by such functions is indeed dynamically invariant under the flow, and analyze its evolution.
The Fourier coefficients corresponding to (\ref{ansatz-gg}) are
\beq
\alpha_0=b\quad\mbox{and}\quad \alpha_n=(a+bp) p^{n-1} \quad\mbox{for}\,\,\, n\geq 1.
\label{singlepole}
\eeq
Expressed through the parameters of our ansatz,
the mass and momentum take the form
\begin{equation}\label{charges}
M=|b|^2+\frac{|a+bp|^2}{1-|p|^2},\quad P=\frac{|a+bp|^2}{(1-|p|^2)^2}\,.
\end{equation}
Substitution of \eqref{ansatz-gg} into (\ref{GGgen}) produces exactly three distinct dependences on $z$ on both sides, namely, $(1-pz)^{-2}$, $(1-pz)^{-1}$ and a $z$-independent term, resulting in the following three equations for $a(t)$, $b(t)$, and $p(t)$
\begin{eqnarray}
 i \dot a &=& M a,\label{eqGGa}\\
 i \dot b &=& (M+P)b+ P a \bar p,\label{eqGGb}\\
 i \dot p &= & M p + a \bar{b} ,\label{eqGGp}
\end{eqnarray}
where we used the conservations laws \eqref{charges} to simplify some terms.
These equations describe a three-dimensional Liouville-integrable system which can be easily solved (in \cite{GG2} an explicit formula for the general solution  was derived using harmonic analysis tools). One first integrates (\ref{eqGGa}) to find $a$, whereupon the remaining equations become linear. The solution corresponding to initial conditions $a(0)=a_0, b(0)=b_0, p(0)=p_0$ (for simplicity assumed  to be real) is
\begin{eqnarray}
 a(t) &=& a_0 e^{-i M t}, \label{adotGG}\\
 b(t) &=& \left(b_0 \cos{\omega t} - i\, \frac{b_0(M+P)+ 2 a_0 p_0 P}{2\omega}\, \sin{\omega t} \right)\, e^{-\frac{i}{2}(M+P) t}, \label{bdotGG}\\
 p(t) &=& \left(p_0 \cos{\omega t} - i\, \frac{p_0 (M+P) + 2 a_0 b_0}{2\omega}\, \sin{\omega t}\right)\, e^{-\frac{i}{2}(M-P) t},\label{pdotGG}
\end{eqnarray}
where $\omega= \frac{1}{2} \sqrt{(M+P)^2-4 P a_0^2}$ and the expressions for $M$ and $P$ in terms of $(a_0,b_0,p_0)$ are given in \eqref{charges}.
Thus, the solution is quasiperiodic for all initial conditions, while
 the corresponding `mass' spectrum,
\beq
|\alpha_0|^2=M-P (1-|p|^2)\quad\mbox{and}\quad |\alpha_n|^2=P (1-|p|^2)^2 |p|^{2(n-1)} \quad\mbox{for}\,\,\, n\geq 1,
\eeq
is expressible through $|p|^2$ and the conserved quantities, and hence exactly periodic.
In the special case $M=P=a_0^2$, we have $\omega=0$ and then
\begin{equation}\label{special}
 a(t)=a_0 e^{-iMt},\quad b(t)=-a_0 p_0 e^{-iMt},\quad p(t)=p_0.
\end{equation}

The evolution of two-mode initial data $a_0=1, b_0=2\varepsilon, p_0=0$ was used in \cite{GG} to illustrate a very interesting instability phenomenon. For these data, \eqref{pdotGG} gives
\begin{equation}\label{2mode_gg}
 p(t)=-\frac{i}{\sqrt{1+\varepsilon^2}} \, \sin(\omega t)\, e^{-2i \varepsilon^2 t},\qquad \omega=2\varepsilon \sqrt{1+\varepsilon^2}\,,
\end{equation}
hence, in the regime of $\varepsilon \rightarrow 0$, $|p(t_n)| \sim 1-\varepsilon^2/2$ for a sequence of times $t_n=\frac{n \pi}{2\omega}$, meaning that the energy goes periodically to arbitrarily high modes. Put differently, although all solutions are quasiperiodic, their radius of analyticity is not uniformly bounded from below. Note that \eqref{2mode_gg} implies instability of the one-mode stationary state $u(t,z)=a_0 e^{-i M t}\,z$.

As mentioned above, one can consider more general meromorphic solutions with an arbitrary number of time-dependent simple poles outside the unit disk
\beq
u(t,z)=\sum\limits_{k=1}^{N} \frac{b_k(t)}{1-p_k(t) z},\quad |p_k|<1,\quad N\geq 2.
\label{multipole}
\eeq
One can convince oneself that the ansatz is consistent in the sense that substituting it in (\ref{GGgen}) results only in $z$-dependences of the form $(1-p_kz)^{-1}$ and $(1-p_kz)^{-2}$ on both sides, producing $2N$ ordinary differential equations for the $2N$ functions $b_k$ and $p_k$. It is possible to consistently impose $p_N=0$, producing an ansatz with one parameter less. The single pole ansatz (\ref{ansatz-gg}) is precisely of such a form.
The equations of motion within the subspace defined by (\ref{multipole}) were written down explicitly in \cite{GG}.
As in the simplest case of the three-dimensional invariant subspace discussed above, the dynamics on all finite-dimensional invariant subspaces is Liouville-integrable and bounded, and hence quasiperiodic.

 Among the quasiperiodic solutions there exist special solutions with time-independent amplitudes $|\alpha_n|$. Such solutions, which we call stationary states, were classified in \cite{GG} in the case of finite dimensional invariant manifolds. They are given by either finite Blaschke products
\begin{equation}\label{blaschke-GG}
 u(t,z)=c\,  e^{-i |c|^2 t}\,\prod\limits_{k=1}^N \frac{\bar p_k -z}{1-p_k z}, \quad |p_k|<1,\quad c\in \mathbb{C},
\end{equation}
 or
\begin{equation}\label{N-GG}
 u(t,z)= \frac{c z^{\ell}}{1-p^N z^N} \, e^{-i \lambda t},\quad p(t)=p(0) e^{-i\omega t},\quad \lambda=\frac{|c|^2}{(1-|p|^{2N})^2},\quad N\omega=\frac{|c|^2}{1-|p|^{2N}},
\end{equation}
where $N\geq 1$ and $\ell\leq N-1$ are nonegative integers.

It was recently shown in \cite{GG3} that outside the finite-dimensional invariant subspaces there exists a dense set of smooth solutions whose radius of analyticity tends to zero for a sequence of times $t_n\rightarrow \infty$. Consequently, the Sobolev norms
with $s>1/2$ are unbounded.
Such a weakly turbulent behavior is somewhat surprising in a completely integrable model. The coexistence of integrability and turbulence is possible because the infinitely many conserved quantities are too weak to control higher regularity properties of solutions.


\section{Three-dimensional invariant subspaces of the conformal flow}
We now return to the conformal flow \eqref{flow}
and attempt to treat it along the lines applied above to the cubic Szeg\H{o} equation. While the treatment is closely parallel, we shall be able to immediately recover only a part of the results available for the cubic Szeg\H{o} equation. Using the generating function $u(t,z)$ and its conjugate $\tilde u(t,z)$ defined as in (\ref{defu}) and (\ref{defutilde}), we find that the system \eqref{flow} is equivalent to the integro-differential equation
\beq
i\del_t\del_z(z u)=\frac1{2\pi i}\oint\limits_{|w|=1} \frac{dw}{w}\tilde u(w) \left(\frac{w u(w)-z u(z)}{w-z}\right)^2.
\label{confgen}
\eeq
The following summation formula has been used while obtaining this complex plane representation:
\beq
\sum_{j=0}^\infty \sum_{k=0}^{n+j} [\min(n,j,k,n+j-k)+1]\rho^j \theta^k=\frac{1-\te^{n+1}}{(1-\te)(1-\rho)(1-\te\rho)}=\frac{1+\te+\cdots+\te^n}{(1-\rho)(1-\te\rho)}.
\label{mastersum}
\eeq
Derivations of \eqref{confgen} and \eqref{mastersum} and further comments are given in appendix \ref{smpl}.

Although equation (\ref{confgen}) looks like a somewhat more elaborate version of (\ref{GGgen}), we do not see an immediate way to characterize its meromorphic solutions and defer it to future work. Nonetheless, the lowest-dimensional nontrivial invariant subspace is easy to construct, as we shall explicitly demonstrate now.

The ansatz relevant for the three-dimensional invariant subspace, analogous to \eqref{singlepole}, is
\beq
\alpha_n=(b+a n)p^n,\quad
\label{confansatz}
\eeq
where $b, a, p$ are complex-valued functions of time. The corresponding generating function has the following combination of poles outside the unit disk
\beq
u(t,z)=\frac{b(t)-a(t)}{1-p(t)z}+\frac{a(t)}{(1-p(t)z)^2}=\frac{b(t)+(a(t)-b(t)) p(t) z}{(1-p(t)z)^2}.
\label{confansatzu}
\eeq
We note that the two-mode initial data ($\alpha_n=0$ for $n\ge 2$) are accommodated within this ansatz as a special limiting case $p\to 0$ with $b$ and $(a-b)p$ finite. All of our general statements about the dynamics on the three-dimensional invariant subspace apply to solutions starting with such initial configurations.

While we could have used the complex plane representation (\ref{confgen}) to establish the validity of our ansatz, it is instructive to apply brute force summations in this case. Substitution of (\ref{confansatz}) into (\ref{flow}) yields
\begin{align}
&i(1+n)\left(\dot b +\dot a n +n(b+a n)\frac{\dot p}{p}\right)\label{harmonic_eq}\\
&= \sum_{j=0}^\infty \sum_{k=0}^{n+j} [\minPlusOne{n,j,k,n+j-k}](\bar b+\bar a j)(b+a k)(b+a(n+j-k))|p|^{2j}.\nonumber
\end{align}
Note that the $p^n$ factor has consistently cancelled on the two sides. It remains to show that different $n$-dependences on both sides can be matched and produce a sufficiently small number of equations.
All summations in (\ref{harmonic_eq}) can be performed by applying the relation
\beq
\sum_{j=0}^\infty \sum_{k=0}^{n+j} [\minPlusOne{n,j,k,n+j-k}] j^K k^L |p|^{2j}= \left(\rho\del_\rho\right)^K\left(\te\del_\te\right)^L\frac{1+\te+\cdots+\te^n}{(1-\rho)(1-\te\rho)}\Bigg|_{\hspace{-1.7mm}\begin{array}{l}\ssty\te=1\vspace{-2.5mm}\\ \ssty\rho=|p|^2\end{array}}\hspace{-3mm},
\label{masterdiff}
\eeq
which follows from (\ref{mastersum}). While we give explicit expressions for the sums involved in (\ref{harmonic_eq}) in appendix \ref{summatns}, the only thing one needs to know about these sums to establish the validity of our ansatz is that they are all polynomials of degree $L$ in $n$, times $(n+1)$, the latter factor coming from Faulhaber's sums $1+2^l+\cdots+n^l$, with $l\le L$, originating from $1+\te+\cdots+\te^n$ in (\ref{masterdiff}). Therefore, explicit counting tells us that, upon substituting the summation formulas of appendix \ref{summatns} in (\ref{harmonic_eq}), the factor $(n+1)$ will cancel on the two sides, leaving behind a statement that two quadratic polynomials in $n$ equal each other. Matching the coefficients of these polynomials produces three ordinary differential equations for three functions $b$, $a$, $p$, confirming the validity of our ansatz.

The explicit equations for $b$, $a$ and $p$ are given by
\begin{align}
&\frac{i\dot p}{(1+y)^2}=\frac{p}{6}\left(2y |a|^2 + \bar b a\right),\label{eqp}\\
&\frac{i\dot a}{(1+y)^2}=\frac{a}{6}\left(5|b|^2+(18y^2+4y) |a|^2+ (6y-1)\bar b a+10y \bar a b  \right),\label{eqr2}\\
&\frac{i\dot b}{(1+y)^2}=b \left(|b|^2+(6y^2+2y) |a|^2 +2y b \bar a\right)+
a \left(2y |b|^2+(4y+2)\s^2 |a|^2 +y^2 \bar b a\right), \label{eqr1}
\end{align}
where we have introduced the following notation, which will turn out to be useful later:
\beq
y=\frac{|p|^2}{1-|p|^2}.
\label{sigmadef}
\eeq
The conservation laws \eqref{charge} and \eqref{lenergy} take the form
\begin{align}
Q=&(1+y)^2\left(|b|^2+4y \Re(\bar b a)+2\s(3y+1)|a|^2\right),\label{Q3d}\\
E=&(1+y)^2\left((1+2y)|b|^2+4y(3y+2)\Re(\bar b a)+ 4y(6y^2+6y+1) |a|^2\right).\label{EQ3d}
\end{align}

We remark that if $u(t,z)$ is a solution of the conformal flow, so is $z^{N}u(t,z^{N+1})$ for any nonegative integer $N$. This automatically generates an infinite number of other three-dimensional invariant subspaces of the conformal flow, in which only subsets of modes are activated.

We shall return to equations (\ref{eqp}-\ref{eqr1}) in section \ref{dynamical}, and analyze the dynamics within the three-dimensional dynamically invariant subspace (\ref{confansatz}) explicitly. Before proceeding in that direction, we shall discuss special solutions for which $|\alpha_n|$ are time-independent.


\section{Stationary states}

The conformal flow (\ref{flow}) admits solutions of the form
 \begin{equation}\label{stac}
  \alpha_n(t)=A_n e^{-i \lambda_n t},
 \end{equation}
 where the frequencies $\lambda_n$ and complex amplitudes $A_n$ are time-independent, and $\lambda_n$ are linear in $n$, that is $\lambda_n=\lambda - n \omega$ for some real $\lambda$ and $\omega$. For such solutions there is no energy transfer between the modes, hence we call them {\em stationary states}. The cubic Szeg\H{o} equation admits a variety of stationary states \cite{GG}; such solutions have also been considered in the context of AdS stability \cite{BBGLL,GMLL,CEJV}, where they were referred to as `quasiperiodic' solutions.

 Substituting \eqref{stac} into \eqref{flow} we get a nonlinear `eigenvalue' problem
\begin{equation}\label{sys}
 (n+1) (\lambda-n \omega) A_n=\sum\limits_{j=0}^{\infty} \sum\limits_{k=0}^{n+j} [\minPlusOne{n,j,k,n+j-k}] \,\bar A_j A_k A_{n+j-k}\,.
\end{equation}
The simplest solutions of this algebraic system, easily seen by inspection, $A_n=c\, \delta_{Nn}$, $\lambda=|c|^2$, $\omega=0$ for $c\in \mathbb{C}$, give
 the one-mode stationary states for any non-negative integer $N$
 \begin{equation}\label{1mode}
  \alpha_n= c\, \delta_{Nn}\, e^{-i |c|^2 t}.
 \end{equation}
We will see shortly that these trivial solutions are the endpoints of two-parameter families of stationary solutions. We note in passing that they are seeds for time-periodic solutions of the original equation \eqref{eq} whose construction will be described elsewhere along the lines of \cite{MR}.

Within the ansatz \eqref{confansatz},
 stationary states take the form
\begin{equation}\label{stac-ansatz}
 b(t)=b(0) e^{-i\lambda t},\quad a(t)=a(0) e^{-i\lambda t},\quad p(t)=p(0) e^{i\omega t}\,.
\end{equation}
Plugging this into the system (\ref{eqp}-\ref{eqr1}), we obtain a system of algebraic equations that can be solved explicitly. The case $a(0)=0$ yields
\begin{equation}
 \label{r1}
 b(t) = c\exp\left(-\frac{i |c|^{2} t}{(1-|p|^{2})^{2}}\right).
\end{equation}
For this solution
$Q=\frac{|c|^2}{(1-|p|^{2})^{2}}$, hence $\lambda=Q$.
The corresponding generating function reads
\begin{equation}
 \label{gen-ground}
 u(t,z) = \frac{c}{1-pz}\,e^{-i\lambda t}.
 \end{equation}

 For nonzero $a$ we get the following two-parameter family of stationary states with $\omega=0$
\begin{equation}
 \label{stac0}
 b(0) = - 2 c |p|^2, \quad a(0)=c (1-|p|^2),\quad \lambda = \frac{|c|^2 |p|^{2}}{(1-|p|^{2})^{2}}
\end{equation}
For this solution
\begin{equation}
 \label{stac0-charges}
 Q = \frac{2 |c|^2 |p|^{2}}{(1-|p|^{2})^2}, \qquad
  E = \frac{4 |c|^2 |p|^{2}(1+|p|^{2})}{(1-|p|^{2})^{3}}, \quad
\end{equation}
hence $\lambda=\frac{1}{2} Q$. The corresponding generating function reads
\begin{equation}
 \label{gen-stac0}
 u(t,z) = c\, \frac{-2 |p|^2 +
  \left(1+\left|p\right|^{2}\right)p z}{(1-pz)^{2}}e^{-i\lambda t}.
\end{equation}
 \vskip 0.2cm
 In addition, for $|p|\leq p_*:=2-\sqrt{3}\approx 0.268$, there is a pair of two-parameter families of stationary states with nonzero $\omega$. The range of $|p|$ is restricted by the condition that $\kappa:=\sqrt{|p|^4-14 |p|^2+1}$ should be real, so we have a tangential bifurcation at $p_*$. These two solutions (labelled by $\pm$) are given by
 \begin{eqnarray}
 b_{\pm}(0) &=& - c(1+5 |p|^{2} \pm \kappa),\quad a_{\pm}(0) = 2 c (1-|p|^2),\\
 \omega_{\pm} & = &\frac{|c|^2}{3} \, \frac{1+|p|^2\pm \kappa}{1-|p|^2},\\
 \lambda_{\pm} & =& \frac{2 |c|^2}{3}\, \left(\frac{3-4|p|^2}{1-|p|^2}
\pm \frac{(3+4|p|^2)\kappa}{(1-|p|^2)^2}\right)\,,
\end{eqnarray}
and their conserved quantities are
\begin{equation}
 Q_{\pm}=\frac{6}{7} (\lambda_{\pm} + \omega_{\pm}),\qquad E_{\pm}=6 \omega_{\pm}\,.
\end{equation}
The corresponding generating function reads
\begin{equation}\label{ge-stacpm}
 u(t,z)=c\, \frac{-(1+5 |p|^2 \pm \kappa) +(3+3|p|^2\pm \kappa) pz}{(1-pz)^2} \, e^{-i \lambda_{\pm} t},
\end{equation}
where $p(t)=p(0) e^{i \omega_{\pm} t}$.

 Outside the three-dimensional invariant subspace given by (\ref{confansatz}) there exist many other stationary states. For instance, we verified that, in close parallel to the cubic Szeg\H{o} equation, any finite Blaschke product
\begin{equation}\label{blaschke}
 u(t,z)=c\, e^{-i |c|^2 t}\, \prod\limits_{k=1}^N \frac{\bar p_k -z}{1-p_k z}
\end{equation}
yields a stationary state.
There are also stationary states where only every $N$th mode is activated, for instance
\begin{equation}\label{N2}
 u(t,z)= \frac{c z^{N-1}}{1-p^N z^N} \, e^{-i \lambda t},\qquad \lambda=\frac{ |c|^2}{(1-|p|^{2N})^2}\,.
\end{equation}
\begin{rem}
For all the above stationary states with $\omega=0$ we have $\lambda=\frac{Q}{N+1}$, where $N$ is the number of zeroes of the generating function (counted with multiplicity). It would be interesting to find a reason that underlies this `quantization' structure.
\end{rem}

It is well known that stationary states of Hamiltonian systems admit a variational characterization. In our case, it follows from the Hamilton equations \eqref{ham-eq} that stationary states \eqref{stac} are the critical points of the functional
\begin{equation}\label{K}
 K:=\frac{1}{2} H - \lambda Q + \omega (E-Q)\,,
\end{equation}
hence stationary states with $\omega=0$ are extrema of $H$ for fixed $Q$, while stationary states with nonzero $\omega$ are extrema of $H$ for fixed $Q$ and $E$. This fact is very helpful in determining stability properties of stationary states, as will be described elsewhere.


\section{Dynamics on the three-dimensional invariant subspace}\label{dynamical}

We now return to equations (\ref{eqp}-\ref{eqr1}) describing the dynamics on the three-dimensional invariant subspace of the conformal flow and demonstrate how to solve them explicitly. Once again, the situation is closely parallel to the three-dimensional invariant subspace of the cubic Szeg\H{o} equation with its underlying one-dimensional periodic motion.

From \eqref{eqp}, we obtain for $y$ defined by (\ref{sigmadef})
\beq
\dot y= \frac{1}{3}y(1+y)^3\Im(\bar b a)\label{eqsigma}.
\eeq
This equation and \eqref{eqr2} imply that
\beq
S=|a|^2y(1+y)^3
\label{S3d}
\eeq
is conserved. This quantity is related to the conformal flow Hamiltonian (\ref{Hconf}) by
\beq
H=Q^2-2S^2.
\label{HS}
\eeq
Equations (\ref{Q3d}), (\ref{EQ3d}) and (\ref{S3d}) can be resolved to express $|b|^2$, $|a|^2$ and $\Re(\bar ba)$ through $Q$, $E$, $S$ and $\s$ as follows
\begin{align}\label{ESQ}
& |b|^2=\frac{2Q-E+3y(Q+2S)}{(1+y)^3},\qquad |a|^2=\frac{S}{y(1+y)^3},\\
&\Re(\bar b a)=\frac{E-Q-2S-2y(Q+6S)}{4y(1+y)^3}.\nonumber
\end{align}
From these relations, together with (\ref{eqsigma}), one obtains
\beq
\dot y^2=-\frac{Q^2+12S^2}{36}\left(y^2+\left(1-\frac{E(Q+2S)}{Q^2+12S^2}\right)y+\frac{(E-Q-2S)^2}{4(Q^2+12S^2)}\right).
\label{sigmaharmonic}
\eeq
Note that this equation has the algebraic form of energy conservation for an ordinary one-dimensional harmonic oscillator. The solution of \eqref{sigmaharmonic} reads
\begin{equation}\label{sigmasol}
 y(t)=B + A \sin(\Omega t+\psi), \qquad \Omega=\frac{1}{6} \left(Q^2+12 S^2\right)^{1/2},
\end{equation}
where $B$ and $A$ are constants depending on $E$, $Q$, and $S$; $\psi$ is a phase that must be determined from the entire set of initial conditions.
 By equation \eqref{ESQ}, this exactly periodic motion is transferred to $|b|^2$, $|a|^2$ and $\Re(\bar ba)$, and hence to the mode energy spectrum
\beq\label{energyspectrum}
|\alpha_n|^2=|b+na|^2|p|^{2n}=\left(|b|^2+2n \Re(\bar ba)+n^2|a|^2\right)|p|^{2n}.
\eeq
One will thus observe exact returns of the energy spectrum to the initial configurations for {\em all} solutions within the three-dimensional invariant subspace of the conformal flow.

The turning points of the periodic motion of $\s$ described by (\ref{sigmasol}), $\s_{\pm}=B\pm A$, are given by the roots of the quadratic polynomial appearing on the right-hand side of \eqref{sigmaharmonic} (the special case $A =0$ corresponds to stationary states). In terms of the energy spectrum \eqref{energyspectrum}, these turning points provide lower and upper bounds for the inverse and direct cascades of energy, respectively. One of the key questions in this context is how large can $y$ grow starting from some small $\s$. To answer this question, we consider the ratio
\beq
\frac{1+y_+}{1+y_-}=\frac{(1+y_+)^2}{(1+y_-)(1+y_+)}\le \frac{(1+y_++y_-)^2}{1+y_-+y_++y_-y_+}.
\label{ratio1}
\eeq
An advantage of the last representation is that it only contains combinations of $y_+$ and $y_-$ directly expressible through the coefficients of the quadratic polynomial on the right-hand side of (\ref{sigmaharmonic}). One thus gets
\beq
\frac{(1+y_++y_-)^2}{1+y_-+y_++y_-y_+}=\frac{4E^2(Q+2S)^2}{(Q^2+12S^2)(E+Q+2S)^2}\le \frac{4(Q+2S)^2}{Q^2+12S^2}.\label{ratio2}
\eeq
We furthermore notice that
\beq
Q=(1+y)^2|b+2y a|^2+2S\ge 2S.
\eeq
Therefore, dividing the numerator and denominator of (\ref{ratio2}) by $Q^2$ and replacing the numerator by its maximum and denominator by its minimum, we obtain a simple uniform upper bound (which very likely can be tightened with extra work)
\beq
\frac{1+y_+}{1+y_-}\le 16.
\eeq
This proves that the transfer of energy to high frequencies (or, equivalently, the growth of the higher Sobolev norms (\ref{sobolev})) is uniformly bounded. We recall from section~\ref{szegosection} that there is no such bound for the Szeg\H{o} flow on the three-dimensional invariant subspace \eqref{ansatz-gg}. In this sense, the conformal flow is much less `turbulent' than the Szeg\H{o} flow. Heuristically, this is not very surprising in view of the fact that the interaction coefficients in \eqref{flowdenom}
\begin{equation}\label{inter-decay}
 \frac{\minPlusOne{n,j,k,n+j-k}}{n+1}
\end{equation}
decay when the recipient mode number $n$ is large and at least one of the source mode numbers is small, undermining the efficiency of energy transfer from low to high frequencies. This is in contrast to the Szeg\H{o} flow where  all coupling coefficients are  equal to 1.


\section{Summary and open questions}

Starting with a naturally defined geometric PDE (\ref{eq-conf2}) on a three-sphere describing a self-interacting conformally coupled  scalar field, we have considered its effective time-averaged weak field dynamics, arriving at the conformal flow (\ref{flow}). The conformal flow is both structurally similar and displays a number of dynamical parallels to the previously known cubic Szeg\H{o} equation (\ref{GGalpha}). Building on these analogies, we have revealed the dynamics on the three-dimensional invariant subspace of the conformal flow described by the ansatz (\ref{confansatz}). Within this subspace, the dynamics is Liouville-integrable and bounded, and hence quasiperiodic, with exactly periodic energy flows (\ref{energyspectrum}) displaying an alternating sequence of direct and inverse cascades.
Unlike the cubic Szeg\H{o} case, the conformal flow dynamics within the three-dimensional invariant subspace does not display turbulent behaviors (not even the `weak weak turbulence' of \cite{KSS}). Other similar three-dimensional subspaces can be immediately constructed, as follows from simple properties of the conformal flow equation. There are furthermore direct parallels between the structure of stationary states of the Szeg\H{o} and conformal flows (these states are special solutions for which no energy transfer between the modes occurs). Of these, the Blaschke product states (\ref{blaschke}) are particularly intriguing, as they clearly demonstrate that the parallels between the Szeg\H{o} and conformal flows extend beyond the three-dimensional invariant subspaces that have been our main focus here.
\vskip 0.2cm
We conclude with a list of open questions that we leave for future investigations:
\begin{itemize}
\item Is the conformal flow (\ref{flow}) Lax-integrable?
\item Is there a way to re-express the conformal flow through a projector similar to the Szeg\H{o} projector? Apart from the conceptual importance, such a projector should make  calculations  easier. The complex plane representation of the conformal flow (\ref{confgen}) can be seen as a first step in this direction.
\item Are there higher dimensional invariant subspaces? The three-dimensional subspace we have considered, with its double-pole generating function (\ref{confansatzu}), and parallels to the cubic Szeg\H{o} case strongly suggest that more elaborate meromorphic ans\"atze should work; however, we have not been able to find them.
\item  Are there weakly turbulent solutions? We have explicitly demonstrated that no turbulent behavior occurs within the three-dimensional invariant subspace (\ref{confansatz}). This does not in principle exclude turbulence for other initial conditions; however, we find it unlikely because the conformal flow generally appears to display less energy spread than the Szeg\H{o} flow.
\item Complete classification of stationary states and analysis of their stability properties would contribute to the overall picture of the conformal flow dynamics.
\item To what extent can the properties of the conformal flow, both explicitly demonstrated and putative, be structurally stable with respect to variations of the flow equation? Do these dynamical patterns have applications to the weak field limit of other related systems? Can the Szeg\H{o} and conformal flows be just two representative members of a large hierarchy of equations?
\item To what extent can  the highly structured dynamics of the conformal flow be transferred to the original conformally invariant cubic wave equation \eqref{eq} on the three-sphere? Standard results on time averaging guarantee that in the  weak field regime the conformal flow accurately approximates our original wave equation for a long but limited time. Can our findings have further implications on longer time scales?

\end{itemize}


\section{Acknowledgments}

P.B.\ gratefully acknowledges the stimulating atmosphere at the trimestre ``Ondes Non Lin\'e\-aires'' in IHES where part of this work was done; he is particularly indebted to Patrick G\'erard for helpful discussions.
The work of P.B., D.H., and M.M.\ was supported by
the Polish National Science Centre grant no. DEC-2012/06/A/ST2/00397. P.B.\ and M.M.\ gratefully acknowledge the support of the Alexander von Humboldt Foundation.  D.H. acknowledges a scholarship of Marian Smoluchowski Research Consortium Matter Energy Future from KNOW funding.

The work of B.C.\ and V.L.\ is supported in part by the Belgian Federal Science Policy Office through the Interuniversity Attraction Pole P7/37, by FWO-Vlaanderen through projects G020714N and G044016N, and by Vrije Universiteit Brussel (VUB) through the Strategic Research Program ``High-Energy Physics''. V.L.\ is supported by a PhD fellowship from the VUB Research Council. B.C.\ thanks the organizers of the NumHol2016 meeting at Santiago de Compostela and of the Nordita program ``Black Holes and Emergent Spacetime'' for hospitality while this work was in progress.

The work of O.E.\ is funded under CUniverse
research promotion project by Chulalongkorn University (grant reference CUAASC). O.E.\ would like to thank Latham Boyle, Stephen Green and especially Luis Lehner for stimulating discussions and warm hospitality during his visit to the Perimeter Institute (Waterloo, Canada). The part of this research conducted during that visit was supported by Perimeter Institute for Theoretical Physics. (Research at Perimeter Institute is supported by the
Government of Canada through the Department of Innovation, Science and Economic Development and by the Province of Ontario through the Ministry of Research and Innovation.)


\appendix

\section{Complex plane representation for the conformal flow}\label{smpl}

We first introduce the generating functions $u$ and $\tilde u$ according to (\ref{defu}) and (\ref{defutilde}).
It is true for any contour enclosing the origin but not enclosing the singularities of $u$, and any contour enclosing all singularities of $\tilde u$ that
\beq
\alpha_m=\frac1{2\pi i}\oint \frac{dz}{z} u(z) z^{-m}, \qquad \bar\alpha_m=\frac1{2\pi i}\oint \frac{dz}{z} \tilde u(z) z^{m}.
\eeq
Substituting these expressions to the right-hand side of conformal flow (\ref{flow}), one gets
\begin{align}\label{geneq}
&i(n+1) \dot\alpha_n =
\frac1{(2\pi i)^3}\oint \frac{ds}{s}\frac{dv}{v}\frac{dw}{w^{n+1}} \tilde u(s) u(v) u(w) \\ &\hspace{4cm}\times\sum_{j=0}^\infty \sum_{k=0}^{n+j} [\minPlusOne{n,j,k,n+j-k}]
\left(\frac{s}{w}\right)^j\left(\frac{w}{v}\right)^k.\nonumber
\end{align}
The summation can be performed using the following master formula
\beq
\sum_{j=0}^\infty \sum_{k=0}^{n+j} [\minPlusOne{n,j,k,n+j-k}]\rho^j \theta^k=\frac{1-\te^{n+1}}{(1-\te)(1-\rho)(1-\te\rho)}=\frac{1+\te+\cdots+\te^n}{(1-\rho)(1-\te\rho)}.
\label{mastersum2}
\eeq
This summation formula can be derived by many methods, for example, by noticing that
\beq
2\min(n,j,k,n+j-k)=n+j-|k-n|-|k-j|,
\eeq
and performing brute force resummations of geometric series. The best way to understand it, however, is essentially combinatorial in nature. The right-hand side of (\ref{mastersum2}) can be rewritten as
\beq
\left(1+\te+\cdots+\te^n\right)\sum_{m=0}^\infty \rho^m \left(1+\te+\cdots+\te^m\right).
\eeq
The coefficient of $\rho^j\te^k$ in this expression is just the coefficient of $\te^k$ in the product
\beq
\left(1+\te+\cdots+\te^n\right)\left(1+\te+\cdots+\te^j\right).
\eeq
Counting the number of relevant pairings of powers of $\te$ in the two polynomials of the product reproduces the coefficient on the left-hand side of (\ref{mastersum2}), thus verifying the summation formula.

Using (\ref{mastersum2}), one can rewrite (\ref{geneq}) as
\beq
i(n+1) \dot\alpha_n =
\frac1{(2\pi i)^3}\oint \frac{ds}{s}\frac{dv}{v}dw\,\tilde u(s) u(v) u(w) \frac{w^{-(n+1)}-v^{-(n+1)}}{(1-s/v)(1-s/w)(1-w/v)}.
\eeq
 (One must have $|s|<|v|$ and $|s|<|w|$ in order for the sum to converge.)
Multiplying by $z^n$ (with $|z|<|v|$ and $|z|<|w|$) and summing over $n$, one obtains
after elementary simplification
\beq
i\del_t\del_z(zu)=\frac1{2\pi i}\oint \frac{ds}{s}\tilde u(s) \left[\frac1{2\pi i}\oint dv \frac{u(v)}{(1-s/v)(v-z)} \right]^2.
\eeq
Since by construction, the $v$-integration contour must not enclose any singularities of $u$, the integral inside the square brackets is simply given by a sum of two residues, resulting in the complex plane representation we have quoted in the main text:
\beq
i\del_t\del_z(zu)=\frac1{2\pi i}\oint \frac{ds}{s}\tilde u(s) \left(\frac{su(s)-zu(z)}{s-z}\right)^2.
\label{confgen2}
\eeq

\section{A few summation formulas}\label{summatns}

We assemble below the explicit summation formulas necessary for deriving (\ref{eqp}-\ref{eqr1}) from (\ref{harmonic_eq}). A crucial feature of these expressions in our context is that all of them are proportional to $(n+1)$.
\begin{align}
&\sum_{j=0}^{\infty} \sum_{k=0}^{n+j} [\minPlusOne{n,j,k,n+j-k}] |p|^{2j}=\frac{n+1}{(1-|p|^2)^2},\\
&\sum_{j=0}^{\infty} \sum_{k=0}^{n+j} [\minPlusOne{n,j,k,n+j-k}] j |p|^{2j}=\frac{2(n+1)|p|^2}{(1-|p|^2)^3},\\
&\sum_{j=0}^{\infty} \sum_{k=0}^{n+j} [\minPlusOne{n,j,k,n+j-k}] k |p|^{2j}=\frac{n(n+1)}{2(1-|p|^2)^2}+\frac{(n+1)|p|^2}{(1-|p|^2)^3},\\
&\sum_{j=0}^{\infty} \sum_{k=0}^{n+j} [\minPlusOne{n,j,k,n+j-k}] j^2 |p|^{2j}=\frac{2(n+1)|p|^2}{(1-|p|^2)^3}+\frac{6(n+1)|p|^4}{(1-|p|^2)^4},\\
&\sum_{j=0}^{\infty} \sum_{k=0}^{n+j} [\minPlusOne{n,j,k,n+j-k}] k^2 |p|^{2j}\nonumber\\
&\hspace{2cm}=\frac{n(n+1)(2n+1)}{6(1-|p|^2)^2}+\frac{(n+1)^2|p|^2}{(1-|p|^2)^3}+\frac{2(n+1)|p|^4}{(1-|p|^2)^4},\\
&\sum_{j=0}^{\infty} \sum_{k=0}^{n+j} [\minPlusOne{n,j,k,n+j-k}] j k |p|^{2j}=\frac{(n+1)^2|p|^2}{(1-|p|^2)^3}+\frac{3(n+1)|p|^4}{(1-|p|^2)^4},\\
&\sum_{j=0}^{\infty} \sum_{k=0}^{n+j} [\minPlusOne{n,j,k,n+j-k}] j^2 k |p|^{2j}\nonumber\\
&\hspace{2cm}=\frac{(n+1)^2|p|^2}{(1-|p|^2)^3}+\frac{3(n+1)(n+3)|p|^4}{(1-|p|^2)^4}+\frac{12(n+1)|p|^6}{(1-|p|^2)^5},\\
&\sum_{j=0}^{\infty} \sum_{k=0}^{n+j} [\minPlusOne{n,j,k,n+j-k}] j k^2 |p|^{2j}\nonumber\\
&\hspace{2cm}=\frac{(n+1)(2n^2+4n+3)|p|^2}{3(1-|p|^2)^3}+\frac{(n+1)(3n+7)|p|^4}{(1-|p|^2)^4}+\frac{8(n+1)|p|^6}{(1-|p|^2)^5}.
\end{align}


\end{document}